\input amssym
\bf \centerline{SUBFACTORS AND CONNES FUSION FOR TWISTED LOOP GROUPS}
\vskip .1in
\centerline{Antony Wassermann, Institut de Math\'ematiques de Luminy}
\vskip .2in
\noindent \bf INTRODUCTION. 
\rm The study of subfactors arising from positive energy representations of loop groups was 
initiated with Vaughan Jones in 1989. In the 1990's I related this to the structures arising from conformal field theory,
in particular vertex algebras, their representations and intertwining operators, so-called ``primary fields''. The link
with subfactors was achieved by the intermediary tool of Connes fusion, a relative tensor product operation on bimodules
over von Neumann algebras. With this tool, the fusion of positive energy representations of $LSU(N)$ at a fixed level 
could be defined analytically and explicitly computed, giving the representations the structure of a modular braided 
tensor category. The computation revolved around the four-point function of certain primary fields and its monodromy
properties, encoded in the Knizhnik-Zamolodchikov differential equation. The computation was initially made using 
primary fields corresponding to arbitrary irreducible representation. In this case the corresponding smeared operators 
yield unbounded operators, which gives rise to certain analytic difficulties. This theory is explained in my 1994 ICM 
article and in Jones' s\'eminaire Bourbaki on this subject from 1995. Subsequently, using the fact that the smeared 
primary fields for the vector representation and its dual were actually bounded, I found a method for treating fusion
which used only these bounded fields, published in 1998. All these computations have served as the starting point for all subsequent  
research in the area between operator algebras and conformal field theory. This approach has been successfully applied by 
my Ph.D. students to other groups and Lie algebras arising in conformal field theory, including the diffeomorphism 
group of the circle/Virasoro algebra (Terence Loke), the loop group of the spin group 
(Valerio Toledano--Laredo) and most recently Neveu--Schwarz algebra, one of the super--Virasoro algebras 
(S\'ebastien Palcoux). 

Another former student Robert Verrill
attempted to treat the case of the twisted loop group of $SU(2N)$, which on the surface seemed 
similar to these other cases but raises problems unresolved in his thesis. The same problem occurs when considering 
another ``twisted'' theory, that of the Ramond algebra, the other super--Virasoro algebra. In the language of conformal 
field theory, in both these cases the vertex algebra defined by the untwisted theory admits an automorphism of period two.
In the loop group case, it corresponds to an outer automorphism $\tau$ of $SU(2N)$ 
coming from a symmetry of the Dynkin diagram
which up to inner automorphisms is simply complex conjugation; for the Neveu--Schwarz algebra, it is just the 
${\Bbb Z}_2$--grading of theory. Igor Frenkel and his students, Haisheng Li and Xiaoping Xu, have developed a theory 
of twisted representations (or modules) of vertex algebras and their intertwiners, perfectly adapted to this situation.
Aspects of the theory had been discussed prior to this by Louise Dolan, Peter Goddard and Paul Montague, in connection 
with the construction of the Monster vertex algebra starting from the period three ``triality'' outer automorphism of
$SO(8)$. Based on this algebraic theory, Dan Freed, Michael Hopkins and Constantin Teleman have extended 
their K--theoretic explanation of fusion rules for loop groups representations to the twisted case, as is to be expected because of
the appearance of twisted loop groups in boundary conformal field theory and open string theory (D--branes). In
this case the outer automorphism defining the twisting can be used to modify the conjugation action of the
underlying group $G$ on itself and the data is encoded in the twisted equivariant K--theory of this space, together 
with its structure as a module over the representation ring $R(G)$ of $G$.  

On the other hand in the 1990s Adrian Ocneanu had developed a theory of modules over the modular 
tensor categories associated with subfactors. Using his combinatorial machinery, Ocneanu determined explicitly 
the possible modules for the tensor categories corresponding to $SU(N)$ for $N=2,3,4$. These categories are the same as
those originally discovered by Jones and Wenzl, coming from Temperley--Lieb algebra and more generally Hecke algebras 
at roots of unity; these appear as the centralizer algebras of the quantum groups of $SU(N)$ at roots of unity, 
the counterpart of loop groups in the theory of exactly solved models in statistical mechanics. Various authors have 
given an abstract formulation of Ocneanu's ideas, in particular Vladimir Turaev and A. Kirillov, 
although without any specific example in mind.

In this announcement we will explain how the positive energy representations of the twisted loop group $L^\tau SU(2N)$ 
provide examples of this general theory. In this case the period two automorphism $\tau$ is most 
easily constructed as follows. Starting from ${\Bbb H}^N = {\Bbb R}^{4N}$ with its standard real inner product, we can 
identify $U(2N)$ with the subgroup of $SO(4N)$ commuting with $\rho({\bf i})$, right multiplication by ${\bf i}$. This 
has a period two automorphism $\tau$ given conjugation by $\rho({\bf j})$ and with fixed point subgroup $Sp(N)$, the 
orthogonal symplectic group. We will also $\tau$ to denote conjugation by $\rho({\bf j})$ on $SO(4N)$, an 
inner automorphism. As for loop groups, there is always an irreducible subfactor given by the 
failure of Haag duality and its structure can be unravelled using Connes fusion. The positive energy representations 
of the twisted loop group can be 
fused with those of the untwisted loop group to give new representations of the twisted loop group. Thus the Grothendieck
group generated by the positive energy representations of the twisted loop group form a module over $R_\ell(SU(2N))$,
the Grothendieck ring for level $\ell$ positive energy representations of the untwisted loop group, the so-called 
{\it Verlinde algebra}. The fusion rules with the fundamental representations of $LSU(2N)$ are given by truncations of
the classical tensor product rules of Sheila Sundaram for tensoring representations of $Sp(N)$ by the 
fundamental representations of $SU(N)$ (the exterior powers of the vector representation).  The fusion of two 
positive energy representations of $L^\tau SU(2N)$ gives a representation of 
$LSU(2N)$. As predicted by the general theory, the fusion of the simplest representation with itself, 
corresponding to the trivial representation of $Sp(N)$, is naturally a *--algebra within the category. Its spectral decomposition 
is a truncation of classical decomposition of the compact symmetric space $SU(2N)/Sp(N)$ corresponding to the period 
two isomorphism, given by the classical Cartan--Helgason theorem. The principal results are described explicitly 
in the next section. In unpublished work, Hans Wenzl has given a construction using quantum groups at roots of 
unity which yields an alternative construction of this module over the tensor category, independent of conformal 
field theory. 

The main results are proved by suitably modifying the techniques used for untwisted loop groups. So far we have 
only developed the unbounded approach since the primary fields creating the twisted sectors from the untwisted ones
need not be bounded. (It is possible that there is a more elementary approach using only bounded operators 
could be found.) A technical analytic tool - the ``phase theorem'' - is required to pass from braiding relations
between two smeared primary fields to realtions between bounded intertwiners given by the phases in their 
polar decomposition. Because the absolute values lie in commuting von Neumann algebras, the notorious problems
in manipulating formal commmutation relations between unbounded operators --- described in VIII.5 of 
Reed and Simon --- can be circumvented. Many parts of the theory 
of Verrill, for example determination of vector primary fields between twisted modules, fit into this theory; 
but these computations need to be extended to other fundamental representations. The explicit construction and analytic 
properties of primary fields depends on the fact that the vacuum representation of $L^\tau SO(4N)$ remains 
irreducible when restricted to $L^\tau U(N)$; this is a consequence of Igor Frenkel's twisted version of the 
boson--fermion correspondence. In terms of nets of von Neumann algebras, note that the local algebras defined by 
the twisted loop group define a double cover of the net defined by the loop group, with monodromy defined 
by the period two automorphism $\tau$, regarded as an automorphism of the original net.

\vskip .2in
\noindent \bf 2. STATEMENT OF MAIN RESULTS. \rm Let $G=SU(2N)$ with involution $\tau$ and set $K=G^\tau = Sp(N)$. We 
define the loop group by
$$LG=C^\infty(S^1,G)$$
with the evident action of the rotation group ${\rm Rot}\,S^1\cong {\Bbb T}$. We consider the 
positive energy projective representations $H$ of  
$LG \rtimes {\rm Rot}\,S^1$, which are irreducible when restricted to $LG$ are irreducible
and which, when restricted to ${\rm Rot}\,S^1$, have a decomposition
$$H=\bigoplus_{n\ge 0} H(n)$$
where 
$$R_\theta \xi = e^{in\theta} \xi$$
for $\xi\in H(n)$ with ${\rm dim}\, H(n) <\infty$, $H(0)\ne (0)$. 
The projective representation $LG\rightarrow PU(H)$ defines
a central extension of $LG$ by ${\Bbb T}$, obtained by pulling back the central extension $U(H)$ of $PU(H)$. 
These are classified by a positive integer $\ell \ge 1$, called the {\it level}. Fixing a level $\ell$, 
there are only finitely many irreducible positive energy representations at that level and these are classified by
the zero energy space $H(0)$. This is the subspace of $H$ invariant under rotation, so it is invariant under the
constant loops $G \subset LG$. It is an irreducible representation of $SU(2N)$ with signature 
$f_1\ge \cdots \ge f_{2N}$ satisfying the constraint $f_1-f_{2N}\le \ell$ and it determines the 
isomorphism class of $H$. (Recall that signatures that differ by a constant 
give the same irreducible representation and that signatures can be identified with Young diagrams.) 
Let $V_f$ be the irreducible representation of $SU(2N)$ with signature $f$ and 
$H_f$ the positive energy irreducible representation with $H_f(0)=V_f$. In particular $V={\Bbb C}^{2N}$ will 
denote the vector representation of $SU(2N)$ and $V_k=\lambda^k(V)$ the $k$th exterior power. Thus 
$$V_f\otimes V_k = \bigoplus_{f<_k g} V_g,$$
where the signature $g$ is obtained by adding $k$ boxes to the signature $f$ so that no two lie in the same row. 
The rule for Connes fusion is similar
$$H_f\boxtimes H_k = \bigoplus_{f<_k g} H_g,$$
except only $g$ satisfying $g_1-g_{2N}\le \ell$ are allowed on the right hand side. Let $d(H_f)$ be the 
{\it quantum dimension} of $H_f$, uniquely specified by $d(H_0)=1$, $d(H_f)>0$ and
$$d(H_f)d(H_k) = \sum_{f<_k g} d(H_g).$$
There is an element $g_0\in G$ such that 
$$d(H_f)= {\rm Tr}_{V_f}\, (g_0).$$
Moreover the subfactor given by the failure of Haag duality in $H_f$ is irreducible and has Jones index $d(H_f)^2$.  

Likewise the twisted loop group is defined by
$$L^\tau G= \{g\in C^\infty({\Bbb R},G): f(x+2\pi) =\tau f(x)\}.$$
with the (double cover of the) rotation group ${\rm Rot}^\tau(S^1)={\Bbb R}/4\pi {\Bbb Z}$ acting by translation.
Again we consider projective representations ${\cal K}$ of the semidirect product 
$L^\tau G\rtimes {\rm Rot}^\tau(S^1)$ which are irreducible on $L^\tau G$ and which when restricted to
${\rm Rot}^\tau \, S^1$ yield a decomposition

$${\cal K}=\bigoplus_{n\ge 0, n\in {1\over 2}{\Bbb Z}} {\cal K}(n),$$
where $R_\theta \xi =e^{in\theta} \xi$ for $\xi \in {\cal K}(n)$ with ${\rm dim}\, {\cal K}(n) <\infty$ and 
${\cal K}(0)\ne (0)$.
The projective representation defines a central extension, classified by a level $\ell \ge 1$ and the
rotation invariant subspace ${\cal K}(0)$ is an irreducible representation of $K=Sp(N)$. Recall that 
these are classified by signatures $h_1\ge \cdots \ge h_N \ge 0$ and in this case, at level $\ell$, must satisfy
the constraint 
$$h_1 + h_2 \le \ell.$$ 
Let $W_h$ be the irreducible representation of $Sp(N)$ with signature $h$ and ${\cal K}_h$ the positive energy
irreducible representation with ${\cal K}_h(0)=W_h$. The classical tensor product rule
$$V_k \otimes W_h = \bigoplus_{h<_p f,\, g<_q f,\, p+q=k} W_g$$
was determined by Sheila Sundaram (and can be verified directly using the Weyl character formula).
For Connes fusion an analogous applies
$${\cal H}_k \boxtimes {\cal K}_h = \bigoplus_{h<_p f,\, g<_q f,\, p+q=k} {\cal K}_g$$
where only $g$ with $g_1 + g_2\le \ell$ are allowed on the right hand side. Again there is a unique 
quantum dimension $d({\cal K}_h)>0$ such that

$$\sum d({\cal K}_h)^2 = \sum d({\cal H}_f)^2,$$
and
$$d{\cal H}_k) d({\cal K}_h) = \sum_{h<_p f,\, g<_q f,\, p+q=k} d({\cal K}_g).$$
There is a constant $C>0$ and an element $k_0\in K$ such that 
$$d({\cal K}_h)= C \cdot {\rm Tr}_{W_h}\, k_0.$$  
By Frobenius reciprocity for bimodules,
${\cal K}_h \boxtimes {\cal K}_g$ is a sum of finitely many $H_f$'s. In particular
$${\cal K}_0 \boxtimes {\cal K}_0 = \bigoplus_{f_1=f_2, f_3=f_4, \dots} {\cal H}_f.$$
Recall that, according to the Cartan--Helgason theorem, 
the spectral decomposition of $L^2(G/K)$ as a $G$--module under left translation is given by
$$L^2(G/K)= \bigoplus_{f_1=f_2, f_3=f_4, \dots} V_f.$$
The subfactor defined by the failure of Haag duality in ${\cal K}_h$ is irreducible and 
has Jones index $d({\cal K}_h)^2$.

\vskip .2in 
\noindent \bf 3. MODULES OVER TENSOR CATEGORIES. \rm Given a unitary tensor category of bimodules $X$, it is 
natural to look for bimodules $Y$ such that $Y\boxtimes Y^*$ is a finite direct sum of irreducible bimodules
in the original tensor category. In that case one by decomposing all bimodules $X\boxtimes Y$ and forming 
the associated Grothendieck group, one obtains a modules over the fusion ring of the original category. As an 
example, take the tensor category of finite--dimensional unitary representations of a finite group $\Gamma$. 
Fixing a 2--cocyle $\omega(g,h)$ of $\Gamma$ with 

$$\omega(h,k)\omega(g,hk)=\omega(g,h)\omega(gh,k),$$
we can consider projective representations of $\Gamma$ with cocycle $\omega$, maps $\pi:\Gamma\rightarrow U(W)$ 
such that

$$\pi(gh) =\omega(g,h)\pi(g)\pi(h).$$
Note that the $\omega$--representations form a module over the ordinary representations.

In this case we can define a corresponding $\omega$--twisted group algebra ${\Bbb C}_\omega[\Gamma]$ 
on which $\Gamma$ acts by conjugation by automorphisms. As a $\Gamma$--module we have

$${\Bbb C}_\omega[\Gamma]=\bigoplus W^*\otimes W$$
summed over equivalence classes irreducible $\omega$--representations $W$. There is moreover a natural coaction

$$\delta:{\Bbb C}_\omega[\Gamma] \rightarrow {\Bbb C}[\Gamma] \otimes {\Bbb C}_\omega[\Gamma],\, g\mapsto g\otimes g$$ 
where ${\Bbb C}[\Gamma]$ is the usual group algebra of $\Gamma$. This coaction is {\it ergodic} in that 
$\delta(a) = 1\otimes a$ if and only if $a$ is a scalar. Such $\Gamma$--algebras are all of the form
${\Bbb C}_{\omega_0}[\Gamma_0]$, where $\Gamma_0$ is a normal subgroup of $\Gamma$ and $\omega_0$ is 2-cocycle of 
$\Gamma_0$; the action of $\Gamma$ is given by $\alpha_g(n) = \mu(g,n)gng^{-1}$, where the pair ($\omega_0$, $\mu$) 
defines an element of the relative group $\Lambda(\Gamma,\Gamma_0)$ of Jones.

Analogous constructions and definitions can be made whenever one starts from a unitary tensor category, 
for example defined by positive energy representations of a loop group at a fixed level,  reducing 
the classification of module categories to that of algebras in the original category admitting an ergodic coaction of
the Hopf algebra $\bigoplus X^*\boxtimes X$ of the category. However, as Wenzl has noted and as the formula for 
${\cal K}_0\boxtimes {\cal K}_0$ shows, the combinatorial structure of the subfactor for ${\cal K}_0$ in fact 
approaches that of the infinite--index inclusion $R^G \subset R^K$, 
for a {\it minimal action} of $G$ on the hyperfinite factor $R$ (see LMS notes), a so--called 
``group--subgroup sufactor''. The subfactor associated with an irreducible $\omega$--representation $W$ 
of $\Gamma$ is just $R^\Gamma \subset (R\otimes {\rm End}\, V)^\Gamma$. In this sense, 
since $K$ is far from being a normal subgroup of $G$, the twisted loop group subfactors seem to combine features of
both these types of subfactor.
\vskip .2in
\noindent \bf 4. POSITIVE ENERGY REPRESENTATIONS OF L$^\tau$SU(2N). 
\rm Given a complex Hilbert space $H$ and a projection $P$ 
on $H$, the complex canonical anticommutation relations 
$$a(f)a(g)^*+a(g)^*a(f)=(f,g)I,\,\,a(f)a(g)+a(g)a(f)=0,$$
have a unique irreducible representation ${\cal F}_P$ with cyclic vector $\Omega$ satisfying
$$(a(g_m)^*\cdots a(g_1)^*a(f_1) \cdots a(f_n)\Omega,\Omega)=\delta_{n,m} \,{\rm det}\, (Pf_i,g_j).$$
The space ${\cal F}_P$ can be identified with {\it fermionic Fock space} and $\Omega$ is called the {\it vacuum vector}.
If a unitary $u\in U(H)$ is such that $[u,P]$ is a Hilbert--Schmidt operator then there is a unitary $U \in U({\cal F}_P)$,
unique up to a phase $z$ with $|z|=1$, such that 
$$Ua(f)U^*= a(uf).$$
This gives a projective unitary representation of this subgroup of $U(H)$ on ${\cal F}_P$ satisfying this 
``quantization condition''. In particular with $H=L^2(S^1,{\Bbb C}^N)$ and $P$ the projection onto Hardy space (with 
values in ${\Bbb C}^N$,  
the natural unitary representation of $LU(N)\rtimes {\rm Diff}\, S^1$ is quantized and gives a positive energy projective
representation of the subgroup $LSU(N)\rtimes {\rm Rot}\, S^1$ at level $1$. The level $\ell$ irreducible representations
are obtained as summands of the $\ell$--fold tensor product. Taking $N=1$ and $H=L^2(S^1)$, 
the corresponding representation of $L{\Bbb T}$ is {\it irreducible}. This is the standard case of fermion--boson duality,
so--called because a projective representation of the Abelian group $L{\Bbb T}$, the product of 
${\Bbb Z}\rtimes {\Bbb T}$ and countably many copies of ${\Bbb R}\times {\Bbb R}$ can be constructed on an 
infinite tensor product of the corresponding Schr\"odinger representations (see the 2006 Cambridge lecture notes for example). More generally
the representation of $LU(N)$ on ${\cal F}_P$ is irreducible, because that of its subgroup $LU(1)^N$ is.

There is a twisted version of the above construction and a corresponding irreducibility statement, 
which as Igor Frenkel pointed out is somewhat 
easier to prove. First note that the above construction depends only on the underlying real structure of the Hilbert 
space $H$ and the comlex structure $J=i(2P-I)$ that $P$ defines on it. Indeed $c(f)=a(f)+a(f)^*$ satisfies 
$$c(f)c(g)+ c(g) c(f) = 2 \, {\rm Re}\, (f,g)$$
and 
$$a(f)={1\over 2} (c(f) -i c(Jf)).$$
Orthogonal transformations $T$ are quantized if and only if $[T,J]$ is Hilbert--Schmidt. In particular this gives
a projective positive energy representation of $LSO(2N)$, irreducible since it is already irreducible 
on the subgroup $LU(N)$. Representations of twisted loop groups are obtained by twisting the underlying Hilbert space. 
The basic example is the twisted loop group $L^\tau {\Bbb T}$ where $\tau(z)=\overline{z}$. 
The real Hilbert space is taken to be the completion of the continuous functions $f:{\Bbb R}\rightarrow {\Bbb C}$ with 
$f(t+2\pi) = \overline{f(t)}$ and an appropriately defined complex structure $J$. The situation is simpler because 
${\Bbb T}^\tau = \pm 1$. The fermionic picture is described infinitesimally by operators $e_r$ with $r\in {1\over 2}
{\Bbb Z}$ with $e_r^*=e_{-r}$ and 
$$e_r e_s + e_s e_r = 2 \delta_{r+s,0}.$$
The infinitesimal generator $d$ for rotations satisfies $[d,e_s]=s e_s$ and there is a unique irreducible representation 
with cyclic vector $\Omega$ satisfying $d\Omega=0$, $e_0\Omega=0$ and $e_r\Omega=0$ for $r>0$. The projective 
representation of $L^\tau {\Bbb T}$ is given infinitesimally by the bosonic operators $a_r$ ($r\in {1/2} + {\Bbb Z}$)
defined by
$$a_r={1\over 2} \sum_{p+q=r} e_p e_q.$$
These satisfy 
$$[a_r,e_s]=e_{r+s},\, [d,a_r]=ra_r$$ 
from which it follows that
$$[a_r,a_s]=r\delta_{r+s,0}.$$
On the other hand the simple combinatorial identity
$$\prod (1+t^m) = \prod {{1-t^{2m}}\over {1- t^m}} =\prod (1-t^{2m-1})^{-1},$$ 
shows that $\Omega$ is a cyclic vector for the $a_n$'s. In particular $L^\tau {\Bbb T}$ acts irreducibly.  

To contruct the level one positive energy representations of $L^\tau SU(2N)$, note that $\rho(j)$ lies in $SO(4N)$ 
so there is an element $X$ in the Lie algebra of $SO(4N)$ such that 
$$\rho(j)= e^{2\pi X}.$$
Conjugation by $e^{tX}$ then gives an isomorphism between $L^\tau SO(4N)$ and $LSO(4N)$ which takes preserves 
positive energy representations. On the other hand $L^\tau SO(4N)$ contains $L^\tau U(2N)$. We claim that, 
as in the untwisted case, the representation of $L^\tau U(2N)$ is irreducible. In fact it suffices to check this 
on the subgroup $L^\tau U(2)^N$, which reduces to the case $N=1$ where it is an immediate consequence of the twisted
fermion--boson correspondence. On the other hand because 
$L^\tau SU(2N)\cdot L^\tau U(1)$ is an index two subgroup of $ L^\tau U(2N)$, the restriction of ${\cal F}_P$ breaks 
up as the direct sum of two distinct irreducible representations of the subgroup. This gives the explicit construction
of positive energy representations of $L^\tau SU(2N)$ at level one; as Agrebaui showed in 1985, 
level $\ell$ representations all arise by decomposing the $\ell$--fold tensor products of level one representations.  
\vskip .2in
\noindent \bf 5. VERTEX ALGEBRAS, TWISTED MODULES AND INTERTWINERS. \rm In conformal field theory, the mathematical 
language of vertex algebras and their representations is by now fairly well established. As introductions there
are the notes of Peter Goddard on meromorphic conformal field theory and the introductory book of Victor Kac.
Here we give a tailor--made approach which constructs the twisted theory in what seems to be the most 
direct and elementary way possible. Given a particular conformal field theory, for example the positive 
energy representations of a loop group at a fixed level, the vacuume representation $H_0$ can be given the 
structure of a vertex algebra. Thus to each finite energy vector $v\in H_0$, there is ``field'' given by a formal
power series 
$$V_0(a,z)=\sum \varphi(a,n) z^{-n-\delta_a}$$
where $\varphi(a,n):H_0(k) \rightarrow H_0(k-n)$, The field should create the state $a$ from the vacuum:
$$V_0(a,z)\Omega|_{z=0} = a.$$
They satisfy the following ``locality relation''
$$V_0(a,z)V_0(b,w)=V_0(b,w)V_0(a,z)$$
in sense that the matrix coefficients of each side, regarded as formal power series, define 
the same rational function in ${\cal R}={\Bbb C}[z^{\pm 1}, w^{\pm 1}, (z-w)^{-1}]$.   
Moreover
$$[L_{-1},V_0(a,z)]={d\over dz} V_0(a,z), \, \, [L_0,V_0(a,z)]=z{dV_0(a,z)\over dz}  + \delta_a V_0(a,z).$$
The following associativity relation follows automatically:
$$V_0(a,z)V_0(b,w)=V_0(V_0(a,z-w)b,w),$$
again as formal power series developments of the same rational function in ${\cal R}$.
For the loop group $LSU(2N)$, each element $X$ of the Lie algebra of $SU(2N)$ defines elements $X(n) =Xe^{-in\theta}$ 
of the complexified Lie algebra of $LSU(2N)$ and hence a field $X(z)=\sum X(n) z^{-n-1}$. For different $X$ and $Y$ the
fields $X(z)$ and $Y(w)$ commute in the sense of the locality relation above. The associativity relation can be 
used as a way of defining products $X_{(n)}Y$ of fields, making precise the ``operator product expansion'' 
of theoretical physicists; a fundamental result of Dong implies, that if field $X$ and $Y$ are local with respect 
to $Z$, then so is any one of the products $X_{(n)} Y$. This result can be used to set up inductively the correspondence
between fields and states. 

For any other representation ${\cal H}_f$, there are similar operators $V_f(a,z)$ on ${\cal H}_f$ such that
$$V_f(a,z)V_f(b,w)=V_f(b,w)V_f(a,z)=V_f(V_0(a,z-w)b,w).$$
Given a representation ${\cal H}_k$ intertwiners 
$\varphi_{fg}^k(v,z):{\cal H}_g \rightarrow {\cal H}_f $ ($v\in {\cal H}_k$) should satisfy:
$$\varphi_{fg}^k(v,z) V_g(a,w) = V_f(a,w) \varphi_{fg}^k(v,z) = \phi_{fg}^k(V_k(a,z-w)v,z),$$
where the intertwiners can have $\delta$'s which are rational numbers. 

In the case of twisting by a period two automorphism of the loop group and hence the vertex algebra, the definition of
the $V_f(a,z)$'s can still be defined, but the expansion must involve $z^{1/2}$. This is explained extensively by 
Li and Xu. Their approach is purely algebraic so cannot be used to deduce analytic properties of the operators 
$V(a,n)$ or $\varphi(v,n)$. 

There is, however, a direct method of constructing these fields which also yields their analytic properties. The idea is
to construct all possible fields at level one and then take tensor products. The vertex algebra of $LSO(4N)$ at level one
is completely described by the fermionic/bosonic construction. Indeed for rational lattices $\Lambda \subset V$, 
vertex operators can be defined for every state in the bosonic Fock space 
$\ell^2(\Lambda) \otimes S(V)^{\otimes \infty}$. The simplest such operators are ``current'' fields 
$$X(z)=\sum X(n)z^{-n-1}$$ 
corresponding to bosons; the other vertex operators have are essentially Fubini--Venziano vertex operators of the  
$$\varphi_X(z)=z^{X(0)} \exp \int X_-(z) \exp \int X_+(z) $$
where $X_\pm(z)=\sum_{\pm n>0}X(n)z^{-n-1}$. All other operators can be expressed in terms of these two types of operators.
On the other hand the vertex operators $\varphi_X(z)$ have an important factorization property 
$$\varphi_X(z)\otimes \varphi_Y(z) = \varphi_{X\oplus Y}(z).$$
Since when $\|Y\|=1$, $\varphi_Y(z)$ yields a fermionic field, this implies analytic properties of 
general $\varphi_X(z)$. This means that the modes $\phi_X(n)$ map Sobolev spacess of $L_0$ into other Sobolev spaces 
and hence preserve the space of C$^\infty$ vectors for ${\rm Rot}\, S^1$.    

Having constructed the vertex algebra for $LSO(4N)$ at level one, the problem is then explicitly to use this 
to describe the twisted sector and all the intertwiners. However, Li has given a very simple description of how
to pass from untwisted to twisted modules. In fact the field
$$\Delta(z)=z^{X(0)} \exp -\sum_{n>0} X(n) (-z)^n/n$$
satisfies
$$[L_{-1},\Delta(z)]=-{d\Delta\over dz},\,\, \Delta(z)V(a,z)\Delta(z)^{-1} = V(\Delta(z)a,z).$$
The fields $V(\Delta(z)a,z)$ define the fields in the twisted version of the vacuum sector. 
Intertwiners are then defined using the relation of Dolan, Goddard and Montague:
$$V(a,z)b= e^{zL_{-1}} V(b,-z)a.$$
At level one, this gives an explicit construction of all possible fields.
\vskip .2in
\noindent \bf 6. PRIMARY FIELDS AND BRAIDING RELATIONS. \rm Note that intertwiners labelled by untwisted representations
can only go from a untwisted or twisted representation to similar type of representation, while those labelled by 
a twisted representation must switch the type. From the Dolan--Goddard--Montague correspondence, it follows that to determine 
the dimensions of spaces of intertwiners, only those labelled by untwisted representations need be considered. 
From this it follows immediately that, for every twisted representation ${\cal K}_h$, there is a unique intertwiner 
labelled by ${\cal K}_h$ from the vacuum representation ${\cal H}_0$ to ${\cal K}_h$. 
For the representations ${\cal H}_k$, corresponding to the exterior powers $V_k=\lambda^k V $, 
we want to determine the dimension of the space of intertwiners
of charge ${\cal H}_k$ between two twisted representations. As in the untwisted case an intertwiner of charge
${\cal H}_f$ is determined uniquely by $\varphi(v,z)$ for $v\in {\cal H}_f(0)=V_f$. 
These fields are called primary fields. If
$$\varphi(v,z)=\sum_{n\in {1\over 2}{\Bbb Z}}\varphi(v,n)z^{-n-\delta},$$
then the $\varphi(a,z)$ defines a primary field if and only if
$$[X(m),\varphi(v,n)]=\varphi(Xv,m+n),\,\, [d,\varphi(v,n)]=-n\varphi(v,n).$$
Here $X(n)=e^{-in\theta}$ with $n$ in ${\Bbb Z}$ or ${1\over 2}+{\Bbb Z}$ according to whether 
$\tau X=X$ or $\tau X= -X$. As in the untwisted case the primary field is uniquely determined by its
initial term, the restriction of $\varphi(v,0)$ to a map between ${\cal K}_{g}(0)$ and ${\cal K}_{h}(0)$, i.e. an element $T$ in
of ${\rm Hom}_K(V_f\otimes W_g, W_h)$. Thus for the exterior powers, these spaces are at most one--dimensional by
Sundaram's tensor product rule. Her rule in fact gives the exact dimension when $g$ and $h$ are permissible, i.e.
${\cal K}_g$ and ${\cal K}_h$ exist.  For, as for untwisted representations, the condition that $T$ be the initial term 
of a primary field can be expressed as linear conditions given by operators from a copy of ${\goth sl}_2$ 
given the highest/lowest weight vectors for $K$ in the -1 eigenspace of $\tau$ in ${\goth g}$: for the exterior powers
it is easy to check these conditions directly (compare the computations of Beauville and Faltings for the 
Verlinde formula). 

As in the untwisted case, the braiding properties of primary fields follow by considering the reduced four-point function
$$f(z)=\sum (\phi(v_2,n)\phi(v_3,-n)v_4,v_1)z^n$$
for lowest energy vectors $v_i$. As Verrill showed, considered as a function with values in 
$U={\rm Hom}_K (V_2\otimes V_3\otimes V_4,V_1)$, this satisfies the twisted Khnizhnik--Zamolodchikov ordinary 
differential equation:
$$f(z)=\left( {A\over z} + {B\over 1-z} + {C\over (1-z)^{1/2}}\right) f(z)$$
where $A$, $B$ and $C$ are operators on $U$ together with boundary condition $D_i(z) f(z)=0$ 
where the $D_i$ are analytic with values in $U^*$. Moreover reduced 2--point functions of this
form exhaust all such functions. The same is true for products of primary fields with $V_2$ and
$V_3$ reversed and $z$ replaced by $z^{-1}$. As in the untwisted case the monodromy of the ordinary 
differential equation from 0 to $\infty$ yields the braiding relation, which passing to smeared primary fields 
takes the form
$$\phi^g_{g0} (a)\phi^{k^*}_{0k} (b) = \sum \lambda_h \phi^{k^*}_{gh}(e^{\mu_h}\cdot b) \phi^g_{hk}(e^{\nu_h}\cdot a),$$
where $\lambda_h$, $\mu_h$ and $\nu_h$ are constants,  $a$ and $b$ are supported in the upper and lower 
semicircle of the unit circle and $e^\alpha(\theta)=e^{i\alpha\theta}$. 
Similarly we have the Abelian braiding
$$\phi^g_{hk}(a)\phi^k_{k0} b) =\varepsilon_h \phi^k_{h g }(e^{\delta_h} \cdot b)\phi_{g0}^g(e^{-\delta_h}\cdot a),$$
where $\varepsilon_h$ and $\delta_h$ are constants.
\vskip .2in 
\noindent \bf 7. LOCAL TWISTED LOOP GROUPS. \rm The local theory of twisted loop groups 
is particularly easy to deduce from the untwisted theory because of the construction by embedding in $LSO(4N)$. 
For an interval $I\subset S^1$ and $G=SU(2N)$, the local twisted loop group $L_I^\tau G$ is defined as the 
subgroup of loops such that $f(t)=I$ for $t\notin I$. 
\vskip .1in
\noindent \bf Generation by exponentials. \it Every element $g$ of $L^\tau G$ is a product of exponentials in 
$L^\tau{\goth g}$.
\vskip .05in
\noindent \bf Proof. \rm Note that if $h\in L\tau G$ with $\|h - I\|_\infty <1$, then $h=\exp X$ with
$X=\log I-(I-h)$ in $L^\tau{\goth g}$. So it suffices to find $g_0=I$, $g_1$, ... , $g_{n-1}$, $g_n=g$ with 
$g_i$ in $L^\tau G$ and $\|g_i - g_{i-1}\|_\infty<1$ for $1\le i \le n$. For then $g_{i-1}^{-1}g_i= \exp X_i$ and
$g=\exp X_1 \cdots \exp X_n$. Since smooth loops are uniformly dense in continuous loops, 
it suffices to find such a chain with $g_1,\dots,g_{n-1}$ continuous. Now if $g(0)=g(2\pi)$ is fixed by $\tau$, 
take a path $g^\prime$ joining $I$ to $g(0)$ in $K$.
   If $g(0)\ne \tau g(2\pi)$, then we can take $g^\prime$ to be
a path in $G$ from $I$ to $g(0)$. Define a map of the boundary of the 1--simplex into $G$ as $g^\prime$ on 
the first side emanating from a vertex, $\tau g^\prime$ on the second side and $g$ on the third side, 
opposite the vertex. Since $G$ is simply connected, the map can be extended to a continuous 
map of the 1-simplex into $G$. Each segment parallel to the third side yields a continuous map $h_t$ ($0\le t\le 1$) 
with $h_t(0)=h_t(2\pi)\in K$. By construction $h_0=I$ and $h_1=g$. Now take $g_i=h_{i/n}$ for $n$ sufficiently large.

\vskip .1in
\noindent \bf Covering property. \it If $S^1=\bigcap_{k=1}^n I_k$, then $L^\tau G$ is generated by the subgroups
$L_{I_k}^\tau G$. 
\vskip .05in
\noindent \bf Proof. \rm By the previous result, it is enough to show that 
every exponential $\exp X$ lies in the group generated by the $L_{I_k}^\tau G$. Let $\psi_k\in C^\infty(S^1)$ be a 
partition of the identity subordinate to $(I_k)$. Then $X=\sum \psi_k\cdot X$, a sum of commuting elements, so that
$$\exp X = \exp \psi_1 \cdot X \cdots \exp \psi_n \cdot X.$$   
\vskip .1in
\noindent \bf Local equivalence. \it Under the natural identification between $L_IG$ and $L^\tau_IG$, 
the irreducible positive energy projective representations at level $\ell$ become unitarily equivalent, 
generating hyperfinite factors of type $III_1$.  
\vskip.05in
\noindent \bf Proof. \rm The automorphism $\tau$ is the restriction of an inner automorphism ${\rm Ad} x$ of $SO(4N)$ 
where $x=\rho(j)$, also denoted by $\tau$. There is a natural inclusion of $L^\tau G$ in $L^\tau SO(4N)$. 
Choose a Lie algebra element $X$ such that $x=\exp 2\pi X$, let $h(t)=\exp tX$. Then 
the map $g\mapsto hgh^{-1}$ gives an isomorphism between $L^\tau SO(4N)$ and $LSO(4N)$, preserving the 
class of positive energy representations at level $\ell$. This yields an embedding $\psi$ 
of $L^\tau G$ in $LSO(4N)$ and the 
level $\ell$ representations of $L^\tau G$ are obtained by decomposing the $\ell$th tensor power of the 
vacuum representation of $LSO(4N)$. Let $\psi_1$ be the natural embedding of $LG$ in $LSO(4N)$. We know that 
$\psi$ is given by the formula $\psi(g)=hgh^{-1}$. For the interval $I$ we can 
find $h_1\in L_JSO(4N)$ with $J$ a proper interval containing $\overline{I}$ and with $h=h_1$ on $I$. 
It follows that $\psi(g)=h_1 \psi_1(g) h_1^{-1}$ for $g\in L_I^\tau G \equiv L_IG$. The level $\ell$ representations of
$LG$ and $L^\tau G$ arises by decomposing the $\sigma =\pi_0^{\otimes \ell}$ of $SO(4N)$ 
according to these two embeddings. Thus $U=\sigma(h_1)$ gives a unitary intertwiner between the two representations.
Let $M=\sigma(\psi(L_I^\tau G))^{\prime\prime}$ and $M_1=\sigma(\psi_1(L_IG))^{\prime\prime}$. Let $P$ and $P_1$ be 
projections onto irreducible summands $\pi$ of $\sigma\circ \psi$ and $\pi_1$ of $\sigma\circ \psi_1$. 
Thus $P$ lies in $M^\prime$ and 
$P_1$ in $M_1^\prime$. Since $M= UM_1U^*$, it follows that $M_1^\prime =UM^\prime U^*$, so that $UPU^*$ lies in $M_1^\prime$.
Since $M_1$ is a type III factor, $UPU^*$ and $P_1$ are unitarily equivalent in $M_1^\prime$. Hence we can find a unitary
$V\in M_1^\prime$ such that $P_1=VUPU^*V^*$. The unitary $VU$ then implements the unitary equivalence 
between the restrictions $\pi$ and $\pi_1$ to $L_IG\equiv L_I^\tau G$.

\vskip .1in

\noindent \bf von Neumann density. \it If $I_1$ and $I_2$ are touching intervals of the circle whose union 
is a proper interval $I$ and $\pi$ is an irreducible positive energy projective representation of $L^\tau G$ on $H$, 
then $\pi(L^\tau_{I_1}G\cdot L^\tau_{I_2} G)$ is dense in $\pi(L_IG)$ in $PU(H)$ for the strong operator topology.

\vskip .05in
\noindent \bf Proof. \rm By local equivalence this reduces to the untwisted case, where the result is known.

\vskip .1in

\noindent \bf Irreducibility. \it If $\pi$ is an irreducible positive energy representation of $L^\tau G$ on $H$
and $N=\pi(L_I^\tau G)^{\prime\prime}$ and $M=\pi(L_{I^c}^\tau G)^{\prime\prime}$, where $I^c$ is the complementary interval,
then $N^\prime \cap M = {\Bbb C}$.
\vskip .05in
\noindent \bf Proof. \rm On removing points $I$ and $I^c$ can be written as the union of intervals 
$I_1,I_2$ and $I_3,I_4$ in anticlockwise order. It suffices to show that the von Neumann algebra generated by 
$\pi(L_I^\tau G)$ and $\pi(L_{I^c}^\tau G$ is all of $B(H)$. But this is just the von Neumann algebra generated
by $\pi(L_{I_k})$ for $1\le k\le 4$. By two applications of the von Neumann density lemma, the von Neumann algebra 
generated by $I_j$ with $j\ne k$ gives the von Neumann algebra corresponding to $I_k^c$. Two of these form 
an open cover of $S^1$ and so their von Neumann algebras generate $B(H)$.

\vskip .2in
\noindent \bf 8. PHASE THEOREM AND CONNES FUSION. \rm
Using hermiticity, the braiding relations for smeared primary fields can be written in the form
$$a_{g0} b_{k0}^* =\sum \lambda_h b_{hk}^*a_{hk},\,\,a_{hk}b_{k0} =\varepsilon_k b_{hg}a_{g0}\eqno{(1)}$$ 
where $a=a_{g0}$ and $b=b_{k0}$ are called the {\it principal parts}. These two braiding relations imply
the {\it transport formula}  
$$ab^*b =\sum |\lambda_h | b_{hg}^*b_{hg} a = (b^\prime)^*b^\prime a,\eqno{(2)}$$
where 
$$b^\prime=(|\lambda_h|^{1/2} b_{hg}).$$
\vskip .1in
\noindent \bf Phase Theorem. \it Both these sets of
intertwiners can be modified so that both the principal parts are {\it unitary operators} and the 
non--principal parts are bounded. 
\vskip .1in

\rm The proof of this will be explained  below. Granted the theorem, if $x:{\cal H}_0\rightarrow {\cal K}_g$ and
$y:{\cal H}_0 \rightarrow {\cal H}_k$ are arbitrary intertwers, their {\it non--principal} parts are defined by
$$x_{ij}=a_{ij}\pi_f(a^*x),\,\, y_{pq}=b_{pq}\pi_q(b^*y).$$
These still satisfy the braid relations and the transport formula. But then to compute the Connes fusion 
${\cal K}_g\boxtimes {\cal H}_k$, we note that by the transport formula 
$$\|x\otimes y\|^2 =(x^*xy^*y\Omega,\Omega) =\sum_h |\lambda_h| (a^*b_{hg}^*b_{hg} a\Omega,\Omega).$$
Hence 
$$U(x\otimes y)=\oplus_h |\lambda_h|^{1/2} y_{hg}x\Omega$$
defines a unitary of ${\cal K}_g\boxtimes {\cal H}_k$ onto $\oplus_h {\cal K}_h$.
\vskip.05in
\noindent (a) {\it The transport formula and braid relations remain true if the operator $a$ is replaced by its phase.}
Let 
$$b^\prime:{\cal K}_g\rightarrow \oplus {\cal K}_h$$
be given by $b^\prime=(b_{hg})$. Let $x=b^*b$ and $y=(b^\prime)^*b^\prime$. Thus
$$ax=ya,\,xa^*=a^*y.$$
Hence, if $X=(I+x)^{-1}$ and $Y=(I+y)^{-1}$, then we have
$$aX=Ya.$$
Let
$$A=\pmatrix{0& 0\cr a & 0 \cr},\,\, B =\pmatrix{X & 0 \cr 0 & Y\cr}.$$
Thus $AB=BA$ on a core $H^\infty$ for $A$, $0\le B\le I$ and $BH^\infty \subseteq {\cal D}(A^*A)$. 
Passing to the closure of $A$, it follows that $AB\supseteq BA$, 
i.e. $B{\cal D}(A) \subseteq {\cal D}(A)$ and $AB=BA$ on ${\cal D}(A)$. But then, taking adjoints, it follows that
$BA^*\subseteq  A^*B$ and hence that $A^*AB \supseteq BA^*A$, i.e. $B{\cal D}(A^*A) \subseteq {\cal D}(A^*A)$
and $A^*AB = BA^*A$ on ${\cal D}(A^*A)$. (This is just the statement that $a^*a$ and $x$ have the sommuting spectral
projections.) So if
$$C= A^*A=\pmatrix{a^*a & 0 \cr 0 & 0\cr},$$
$(I+C)^{-1}$ commutes with $B$. Similarly $(\varepsilon + C)^{-1}$, and hence its square root $(\varepsilon + C)^{-1/2}$
commutes with $B$.  Now $AB= BA$ on ${\cal D}(A)$, so that
$$A(\varepsilon + C)^{-1/2}  B = AB(\varepsilon + C)^{-1/2}= BA(\varepsilon + C)^{-1/2}$$
on ${\cal D}(A)$. Since $(\varepsilon + C)^{-1/2} A$ tends in the strong operator topology to the phase of $A$, 
it follows that the phase of $A$ commutes with $B$. 

Similarly it follows that $A_p=\pi_p((\varepsilon + a^*a)^{-1/2})$ satisfies $A_p b_{pq} =b_{pq} A_q$. Hence
the intertwiners $a_{ij}\pi_j((\varepsilon + a^*a)^{-1/2})$ satisfy the braiding relations, are bounded and as 
$\varepsilon\rightarrow 0$ have a strong operator limit. This limit also satisifies the braid relations, but
now the principal part is the phase of $a$.   
\vskip .05in
\noindent (b) {\it The transport formula and braid relations remain true of the operator $b$ is replaced by its phase.}
This follows by the similar reasoning to (a), but the proof is slightly easier. Let $A_1=a$, 
$A_2=\oplus |\lambda_h|^{1/2} a_{hk}$, $B_1 = b$ and $B_2=\oplus \varepsilon_h |\lambda_h|^{1/2} b_{hg}$. 
The braiding relations then take the form $A_1B_1^* =B_2^*A_2$ and $A_2B_1=B_2A_1$. By assumption, 
the $A_i$'s are partial isometries. On the other hand it is easy to see that these equations 
are unchanged if the $B_j$'s are replaced by their phases.
\vskip .05in
\noindent (c) {\it The transport formula and braid relations also hold with $a$ and $b$ unitary.}
Each set of intertwiners $(c_{ij})$ can be modified in three steps in three steps, preserving the braiding relations.
Firstly $c_{ij}$ can be replaced by $\sum 2^{-n} \pi_i(g_n) c_{ij} \pi_j(u_n)$with $(g_n)$ a dense subgroup of the 
elevant local loop group and $u_n$ unitaries in $\pi_0(L_IG)^{\prime\prime}$ with $u_iu_j^*=\delta_{ij}I$ 
and $\sum u_i^*u_i=I$. Secondly passing to the phase, we many assume that $cc^*=I$.  Finally, we can replace $c_{ij}$ by
$c_{ij}\pi_j(u)$ where $u$ is a partial isometry in $\pi_0(L_IG)^{\prime\prime}$ with $uu^*=I$ and $u^*u=c^*c$. 
\vskip .1in
\noindent \bf 9. FUSION RULES. \rm In the previous section upper bounds were obtained for fusions of twisted 
representations with the untwisted fundamental representations. In order to show these upper bounds are attained, 
we provide an algebraic model for these fusion rules. (As in the case of the Verlinde formula for fusion of 
untwisted representations, the fusion coefficients can also be expressed in terms of the classical 
tensor product coefficients by incorporating corrections given the action of a suitable affine Weyl group.) 

We start by recalling that the characters
of $SU(2N)$ are labelled by signatures $f_1\ge f_2\ge \ge \cdots \ge f_{2N}$. On a matrix with diagonal entries
$z_i$, the corresponding character is given by
$$\chi_f(z)={{\rm det} z_j^{f_i+2N-i} \over \prod_{i<j} (z_i-z_j)}.$$
Diagonal matrices in $Sp(N)$ have entries $\zeta_1^{\pm 1}, \dots, \zeta_N^{\pm 1}$ and the characters of irreducible 
representations are indexed by signatures $h_1\ge \dots h_N\ge 0$, with
$$\psi_h(\zeta)={{\rm det }\, \zeta_j^{f_i +N -i +1} -\zeta_j^{-f_i-N+i-1}\over \prod(\zeta_i-\zeta_i^{-1}) \prod_{i<j} 
(\zeta_i+\zeta_i^{-1} -\zeta_j-\zeta_j^{-1})}.$$
Now let 
$$S=\{g:g_1\ge g_2\ge \cdots \ge g_N\ge 0,\, g_1\le\ell/2, \,\,\hbox{$g_i$ all integers or half--integers}\}.$$
Let $D(g)\in Sp(N)$ be the diagonal matrix with entries $\zeta_i^{\pm 1}$ with 
$\zeta_i=\exp(2\pi i(g_i + N+1/2-i)/(2N+\ell))$ for $g\in S$. If $h_1+h_2=\ell +1$, then $\psi_h(D_g)=0$ since the first
two columns of the determinant in the numerator vanish. Let ${\cal T}=\{D(g):g\in S\}$ and let 
$\theta:R(Sp(N))\rightarrow \rightarrow {\Bbb C}^{{\cal T}}$ denote the evaluation or restriction map. 
Let ${\cal S}=\theta(R(Sp(N))$ and let $\theta_h=\theta(W_h)$. The $\theta_h$'s with $h$ permissible 
(i.e. $h_1+h_2\le \ell$) are clearly closed
under multiplication by the $\theta(\chi_k)$. Moreover the fundamental representations of $Sp(N)$ are just
$\chi_1$, $\chi_2-\chi_1$, $\chi_N-\chi_{N-1}$ where $\chi_k$ is the character of $\Lambda^k{\Bbb C}^{2N}$. It follows that
the $\theta(\psi_k)$'s generate ${\cal S}$ and hence that the ${\Bbb Z}$--linear span of the $\theta_h$'s with $h\in S$
must equal ${\cal S}$. The characters $\chi_k$ with $1\le k\le N$ distinguish the points of ${\cal T}$. Hence 
${\cal S}_{\Bbb C}$ is a unital subalgebra of ${\Bbb C}^{\cal T}$, so must equal ${\Bbb C}^{\cal T}$. It follows that
the $\theta_h$'s with $h$ permissible form a 
${\Bbb Z}$--basis of ${\cal S}$. We claim that ${\rm ker}\,\theta$ is the ideal
in $R(Sp(N))$ generated by $\theta_h$'s with $h_1+h_1=\ell +1$. Indeed if $I$ is the ideal they generate, then 
$R(Sp(N))/I$ is spanned by the image of the $[W_h]$'s with $h$ permissible, 
from the tensor product rules with the $V_k$'s. But
$I\subseteq {\rm ker}\,\theta$ and the $\theta_h$'s are linearly indendent over ${\Bbb Z}$. Hence the images of the 
$[W_h]$'s give a ${\Bbb Z}$--basis of $R(Sp(N))/I$ and therefore $I={\rm ker}\, \theta$.

So far we know that 
$${\cal H}_k \boxtimes {\cal K}_h \le \bigoplus_{h<_pf,\,g<_qf,\, p+q=k} {\cal K}_g.\eqno{(*)}$$
Taking quantum dimensions it follows that
$$ d({\cal H}_k)d({\cal K}_h) \le \bigoplus_{h<_pf,\,g<_qf,\, p+q=k} d({\cal K}_g).$$
On the other hand $d({\cal H}_k)=\chi_k(D(0))$. The quantities $\psi_h(D(0))$ are also positive and satisfy
$$\chi_k(D(0))\psi_h(D(0))=\sum_{h<_pf,\,g<_qf,\, p+q=k} \psi_g(D(0)).$$
Since the Perron--Frobenius eigenvalue of a matrix with non--negative entries strictly decreases on 
passing to a submatrix, we deduce that
equality occurs in $(*)$ and that $d({\cal K}_h)=C \psi_h(D(0))$ for some $C>0$. It follows immediately that fusion
with more general representations ${\cal H}_f$ can be computed by looking at how the basis elements $\theta_h$ 
of ${\cal S}$ transform under multiplication by $\theta(\chi_f)$. To compute $C$, we note that, from the theory 
of bimodules, 
$$\sum_f d({\cal H}_f)^2 = \sum_h d({\cal K}_h)^2.$$
So 
$$C^2={\sum_f \chi_f(D(0))^2 \over \sum_h \psi_h(D(0))^2}.$$
From the modular transformation rules of the characters of the untwisted and twisted loop group of $SU(2N)$ of 
Kac and Petersen (see Chapter 13 of Kac's book), we have
$$\left(\sum_f \chi_f(D(0))^2\right)^{-1/2} = (2N)^{-1/2}(2N+\ell)^{-N+1/2}2^{N(2N-1)} 
\prod_{1\le i <j \le 2N-1}\sin {(j-i)\pi\over 2N +\ell}$$
and
$$\eqalign{\left(\sum_g \psi_h(D(0))^2 \right)^{-1/2}& =\cr
(2N+\ell)^{-N/2} 2^{N^2}&\prod_{0< i< N, \, i\in 1/2 + {\Bbb Z}}
\sin {2i\pi\over 2N + \ell} \cdot \prod_{0<i<j<N,\,i,j\in 1/2+{\Bbb Z}} \sin {(i+j)\pi\over 2N+\ell} \sin {(j-i)\pi\over
2N+\ell}.\cr}$$
Hence 
$$C^{-1}=(2N)^{1/2} (2N+\ell)^{(N-1)/2}2^{-N(N-1)} {\prod_{1\le 2k+1<2N} \sin^{k+1} {(2k+1)\pi\over 2N+\ell} \prod_{k=1}^{N-1}
\sin^{N-k}{k\pi\over 2N+\ell}
\over \prod_{k=1}^{2N-2} \left(\sin^{2N-k-1} {k\pi\over 2N+\ell}\right)}$$    
\vskip .2in
\noindent \bf REFERENCES. \rm 
\item {1.} A. Beauville, {\it  Conformal blocks, fusion rules and the Verlinde formula.}  
Proceedings of the Hirzebruch 65 Conference on Algebraic Geometry (Ramat Gan, 1993),  
75--96, Israel Math. Conf. Proc., 9, Bar-Ilan Univ., Ramat Gan, 1996.
\item{2.} L. Dolan, P. Goddard, P.S. Montague, {\it Conformal field theory of twisted vertex operators.}  
Nuclear Phys. B  338  (1990),  no. 3, 529--601. 
\item {3.} G. Faltings, {\it  A proof for the Verlinde formula.}  J. Algebraic Geom. 3 (1994),  no. 2, 347--374.
\item{4.} D. Freed, M. Hopkins, C. Teleman, {\it Loop groups and twisted K--theory III}, (2005)
\break http://arxiv.org/abs/math.AT/0312155.
\item{5.} M.R. Gaberdiel, T. Gannon, 
{\it Boundary states for WZW models.}  Nuclear Phys. B  639  (2002),  no. 3, 471--501.
\item{6.} P. Goddard, D. Olive, ``Kac-Moody and Virasoro Algebras: A Reprint Volume for Physicists,'' (1988), 
World Scientific.
\item{7.} P. Goddard, {\it Meromorphic conformal field theory,}
``Infinite-dimensional Lie algebras and groups'' (Luminy, 1988),  
556--587, Adv. Ser. Math. Phys., 7, World Sci. Publ., Teaneck, NJ, 1989.
\item{8.} V.F.R. Jones, {\it Fusion en alg\`bres de von Neumann et groupes de lacets (d'apr\`es A. Wassermann).} (French) 
[Fusion in von Neumann algebras and loop groups (after A. Wassermann)] S\'eminaire Bourbaki, Vol. 1994/95.  
Ast\'erisque  No. 237  (1996), Exp. No. 800, 5, 251--273.
\item{9.} V. Kac, ``Infinite dimensional Lie algebras'', Third edition, Cambridge University Press, 1990.
\item{10.} V. Kac, ``Vertex algebras for beginners'', Second edition, American Mathematical Society, 1998.
\item{11.} A. Kirillov, Jr, {\it On $G$--equivariant modular
categories}, (2004), http://arxiv.org/abs/math/0401119.
\item{12.} H. Li, {\it The theory of physical superselection sectors in terms of vertex operator algebra language,}
\break http://arxiv.org/abs/q-alg/9504026.
\item{13.} H. Li, {\it Local systems of twisted vertex operators, vertex operator superalgebras and twisted modules.}
``Moonshine, the Monster, and related topics'' (South Hadley, MA, 1994),
203--236, Contemp. Math., 193, Amer. Math. Soc., Providence, RI, 1996.
\item{14.} T. Loke, {\it Operator algebras and conformal field theory of the discrete series representations of 
$Diff(S^1)$}, Ph.D. thesis., University of Cambridge, 1994. 
\item{15.} P.S. Montague, {\it Intertwiners in orbifold conformal field theories}, Nuclear Physics B 486 [PM] (1997), 
546--564. 
\item{16.} A. Ocneanu, {\it The classification of subgroups of quantum ${\rm SU}(N)$.}  
``Quantum symmetries in theoretical physics and mathematics'' (Bariloche, 2000),  133--159, 
Contemp. Math., 294, Amer. Math. Soc., Providence, RI, 2002. 
\item{17.} S. Palcoux, ``Unitary discrete series, characters, Connes fusion and subfactors for the Neveu--Schwarz algebra'', Ph.D. thesis, Institut de Math\'ematiques de Luminy, December 2009.
\item{18.} V.B. Petkova, J--B. Zuber, {\it Verlinde NIM-reps for charge conjugate ${\rm sl}(N)$ WZW theory.}  
``Statistical field theories'' (Como, 2001),  161--170, 
NATO Sci. Ser. II Math. Phys. Chem., 73, Kluwer Acad. Publ., Dordrecht, 2002. 
\item{19.} T. Quella, I. Runkel, C. Schweigert, {\it An algorithm for twisted fusion rules.}  
Adv. Theor. Math. Phys.  6  (2002),  no. 2, 197--205.
\item{20.} S. Sundaram, {\it The Cauchy identity for ${\rm Sp}(2n)$}.  J. Combin. Theory Ser. A  53 
(1990),  no. 2, 209--238
\item{21.} S. Sundaram, {\it Tableaux in the representation theory of the classical Lie groups.}  
``Invariant theory and tableaux'' (Minneapolis, MN, 1988),  191--225, IMA Vol. Math. Appl., 19, Springer, New York, 1990.
\item{22.} V. Toledano--Laredo, ``Fusion of Positive Energy Representations of LSpin(2n)'', (1997),  
Ph.D. thesis, University of Cambridge, http://arxiv.org/abs/math/0409044.
\item{23.} V. Turaev, {\it Homotopy field theory in dimension 3 and crossed group--categories}, (2000), 
\break http://arxiv.org/abs/math/0005291.
\item{24.} R.W. Verill, ``Positive energy representations of $L^\sigma$ SU(2r) and orbifold fusion,'' Ph.D. thesis, 
University of Cambridge, July 2001.
\item{25.} A.J. Wassermann, {\it Coactions and Yang-Baxter equations for ergodic actions and subfactors.}  
``Operator algebras and applications'', Vol. 2,  203--236, 
London Math. Soc. Lecture Note Ser., 136, Cambridge Univ. Press, Cambridge, 1988.
\item{26.} A.J. Wassermann, {\it Operator algebras and conformal field theory.}  
Proceedings of the International Congress of Mathematicians, Vol. 1, 2 (Z\"urich, 1994),  
966--979, Birkh\"auser, Basel, 1995.
\item{27.} A.J. Wassermann, {\it Operator algebras and conformal field theory. III. 
Fusion of positive energy representations of ${\rm LSU}(N)$ using bounded operators.}  
Invent. Math.  133  (1998),  no. 3, 467--538.
\item{28.} A.J. Wassermann, {\it Analysis of operators}, Part III lecture notes, Cambridge, Michaelmas 2006. \break http://www.dpmms.cam.ac.uk/~ajw/AO.ps
\item{29.} H. Wenzl, {\it Subgroup type subfactors and twisted Kac--Moody algebras}, Lecture at the Fields Institute, 
Toronto, November 2007. http://www.fields.utoronto.ca/audio/07-08/vonneumann/wenzl/
\item{30.} H. Weyl, ``The classical groups,'' 
The Classical Groups. Their Invariants and Representations. Princeton University Press, Princeton, N.J., 1939.
\item{31.} X. Xu, {\it Intertwining operators for twisted modules of a colored vertex operator superalgebra.}  
J. Algebra  175  (1995),  no. 1, 241--273.
\item{32.} X. Xu, {\it Twisted modules of coloured lattice vertex operators superalgebras.}
Quart. J. Math. Oxford Ser. (2)  47  (1996),  no. 186, 233--259.
\item{33.} X. Xu, ``Introduction to vertex operator superalgebras and their modules.'' 
Mathematics and its Applications, 456. Kluwer Academic Publishers, Dordrecht, 1998. 
\vfill\eject
\end